 \documentstyle[12pt,amstex]{article}

\title{Remarks on central extensions of the Galilei group in $2+1$ dimensions}
\author{Yves Brihaye\\
Dept. of Mathematical Physics,
University of Mons\\
Av. Maistriau, B--7000 Mons, Belgium\\
{\small (e--mail: snuyts2@@vm1.umh.ac.be)}\\
Stefan Giller,$^*$
Cezary Gonera,$^*$
Piotr Kosi\'{n}ski
\thanks{supported by the grant no.\ 458 of the University of
{\L}\'{o}d\'{z}}\\
Dept. of Field Theory,
University of {\L}\'{o}d\'{z}\\
Pomorska 149, 90--236 {\L}\'{o}d\'{z}, Poland\\
{\small (e--mail: pkosinsk@@krysia.uni.lodz.pl)}}
\date{}
\begin{document}
\maketitle
\newpage\vfil
\begin{abstract}
Some properties of central extensions of $2+1$ dimensional Galilei group are
discussed. It is shown that certain families of extensions are isomorfic.
An interpretation of new nontrivial cocycle is offered. Few bibliographical
remarks are included.
\end{abstract}\vfil
\newpage

\section{Introduction}
Some attention has been recently paid to the nonrelativistic symmetry of
$2+1$--dimensional space--time. The relevant Galilei group differs
significantly from its fourdimensiomal counterpart which makes the study
of its mathematical structure quite interesting. Moreover, one can expect
that such an analysis will appear helpful in understanding the propeties
of nonrelativistic systems that are effectivly confined to two spatial
dimensions. In papers \cite{1}, \cite{2} Bose considers the problem of
finding all central extensions of $2+1$--dimensional Galilei group (and
its Lie algebra) and constructs the relevant umitary projective
representations. However, we feel that not all interesting points were
exhausted there. It is the aim of the present short note to add some
futher remarks concerning the central extensions of Galilei group /
algebra in three dimensions. In section II we prove (actualey, the proof is
almost trivial) two theorems indicating that certain families of such
extensions are isomorphic and indicate how they can be used to find the
relevant Casimir operators and to simplify slighty the representation
theory. In section III the explanation is offered for the existence of
additional nontrivial cocycle in three dimensions. Namely, it is shown
that occurence of this cocycle is related to the Thomas precession
phenomenon for threedimensional Lorentz group. Finally, the
telegraphically short section IV contains some bibligraphical notes. This is
because we think the proper credit should be given to authors who obtained
the results contained in Ref. \cite{1}.

\section{\sloppy On central extensions of 2+1-di\-men\-sional
Galilean algebra and
group}
Let $M$, $N_{i}$, $H$ and $P_{i}$ be rotation, boost, time-- and
space-- translation generators, respectively. Bose \cite{1} has proven that
the vector space of central extensions of Lie algebra of Galilei group is
threedimensional. In the notation adopted above the extended algebra reads
\newpage
\begin{equation}
\left.
\begin{array}{ccl}
\left[H,P_{i}\right]&=&0 \\
\left[N_{i},H\right]&=&i P_{i} \\
\left[P_{i},P_{j}\right]&=&0 \\
\left[N_{i},N_{j}\right]&=&i k \varepsilon_{ij}
{\bf 1} \\
\left[M,P_{i}\right]&=&i \varepsilon_{ij} P_{j} \\
\left[N_{i},P_{j}\right]&=&i m \delta_{ij} {\bf 1} \\
\left[M,N_{i}\right]&=&i \varepsilon_{ij} N_{j} \\
\left[M,H\right]&=&i l {\bf 1}
\end{array}
\right\} \label{1}
\end{equation}
the extension being parametrized by three real numbers $m$, $k$, $l$; the
central element has been denoted by {\bf 1}.

Let us denote the above algebra by $g_{kml}$. The following suprisingly
simple result holds:\\
{\bf Theorem I}:\\
{\em Let $m \neq 0$, $l$ -- arbitrary but fixed. Then $g_{kml}$ are
isomorphic, as
Lie algebras, for all k}.\\
{\bf Proof}.\\
Redefine the basis as follows: $X' = X$, $X \neq N_{i}$, $N_{i}' = N_{i} +
\frac{k}{2 m} \varepsilon_{ij} P_{j}$ $\Box$\\
Let us point out that such an isomorphism does not necessarily imply
physical equivalence (cf. Ref. \cite{3}).

As an application we list all Casimir operators for arbitrary
$k$, $m$, $l$. It reads\\
(i) $l = 0$, $m \neq 0$, $k$---arbitrary
\[
\begin{array}{ccl}
C_{1} & = & H - \frac{1}{2 m} \vec{P}^2,\\
C_{2} & = & M - \frac{1}{m} \vec{N} \times \vec{P} - \frac{k}{m} H
\end{array}
\]
(ii) $l$--arbitrary, $m = 0$, $k = 0$\\
\[
\begin{array}{ccl}
C_{1}' & = & \vec{P}^2, \\
C_{2}' & = & \vec{N} \times \vec{P}
\end{array}
\]
(iii) $l$--arbitrary, $m = 0$, $k \neq 0$\\
\[
\begin{array}{c}
C_{1}'' = \vec{P}^2
\end{array}
\]
(iv) $l \neq 0$, $m \neq 0$, $k$--arbitrary -- none.\\
We are not going to give here the detailed proof but rather content
ourselves with few remarks. The case (i) is a straightforward
consequence of Theorem I and the analogy with fourdimensional case;
(ii) and (iii) are easily verified and only (iv) calls for some comments.
Let us put $k = 0$ which, by Theorem I, does not restrict the generality.
Let $C$ be the central element of $g_{0ml}$. It can be written in ``normal``
order as
\begin{equation}
C = \sum_{(\lambda),(\mu), \nu, \rho} c_{(\lambda)(\mu)\nu\rho}
\prod_{i=1}^{2} N_{i}^{\lambda_{i}} \prod_{i=1}^{2} P_{i}^{\mu_{i}}
H^{\nu} M^{\rho}. \label{2}
\end{equation}
Let $\rho_{max}(c)$ be the maximal power of $M$ on the right hand side of
eq.(\ref{2}). Assume that $\rho_{max}(c) > 0$; then
\begin{equation}
C'= C + i (l \rho_{max}(c))^{-1} \cdot \left[C,C_{1}\right] \cdot M = C
\label{3}
\end{equation}
while $\rho_{max}(c') \leq \rho_{max}(c) - 1$; therfore $\rho_{max}(c) =
0$. Now apply the same reasoning with $M$ replaced by $H$ and $C_{1}$ replaced
by $C_{2}$ (with $k=0$) to conclude that $\nu_{max}(c)=0$. We continue this
argument by taking $N_{i}$ and $P_{i}$ instead of $C_{i}$ to show that
$\lambda_{i\, max}(c)=0$ and $\mu_{i\, max}(c)=0$.

Let us now consider the central extensions of Galilei group. The algebras
$g_{km0}$ can be integrated to yield the central extensions $G_{km}$ of
Galilei group $G$ \cite{1} \cite{3}. They can be described as follows. Let
$(\tau, \vec{u}, \vec{v}, R)$ be an element of Galilei group with $\tau$,
$\vec{u}$, $\vec{v}$, $R$ being time translation, space traslation, boost and
rotation, respectively. Then the multiplication rule for $G_{km}$ reads
\[
(\zeta, \tau, \vec{u}, \vec{v}, R) \ast (\zeta', \tau', \vec{u'},
\vec{v'},R') =
\]
\begin{equation}
= (\zeta \zeta'\omega, \tau + \tau', \overrightarrow{R u'} + \vec{v}
\cdot \tau' + \vec{u}, \overrightarrow{R v'} + \vec{v}, R R') \label{4}
\end{equation}
where $\zeta \in {\Bbb C},\, |\zeta| = 1$ and non trivial cocycle is given by
\begin{equation}
\omega=\exp\left(-i m (\frac{\vec{v}^{2}}{2} \tau' + \vec{v} \cdot
\overrightarrow{R u'}
) - \frac{i k}{2} (\vec{v} \times \overrightarrow{R v'})\right) \label{5}
\end{equation}
We adopted here the results of Ref\cite{3}; the corresponding cocycle
differs by coboundary from the one given in Ref\cite{1}.

Theorem I has the following counterpart on the group level.\\
{\bf Theorem II.}\\
{\em Let $m \neq 0$; then all groups $G_{km}$ are isomorphic.}\\
{\bf Proof.}\\
Make the following change of parameters:
\[
u_{i} \rightarrow u_{i} + \frac{k}{2 m} \varepsilon_{ij} v_{j},
\]
the remaining parameters being unaffected. $\Box$\\

Again this result appears to be quite useful. In Ref.\cite{2} the induced
representations of $G_{km}$ have been found following Mackey's method.
However, when attempting to apply this method in straightforward way one
is faced with the following apparent difficulty: there seems to be no
convenient semidirect product structure for $G_{km}$. This difficulty was
overcome in Ref.\cite{2} by considering the extensions of Galilei group $G$
with the help of two central charges and selecting the appropriate
representations. However, in view of our Theorem II it is unnecessary: we
can always assume $k=0$ or $m=0$ and in both cases the semidirect
structure is transparent.
For $l \neq 0$, $g_{kml}$ can be integrated to the central extension $
\tilde{G}_{kml}$ of the universal covering $\tilde{G}$ of Galilei goup.
The relevant group multiplication rule reads
\[
(\zeta, \tau, \vec{u}, \vec{v}, \theta) \ast (\zeta', \tau',\vec{u}',
\vec{v}', \theta') =
\]
\begin{equation}
= (\zeta \zeta' \tilde{\omega}, \tau+\tau',
\overrightarrow{R(\theta) u'}+ \vec{v} \cdot \tau'+ \vec{u},
\overrightarrow{R(\theta) v'}+\vec{v}, \theta+\theta'); \label{6}
\end{equation}
here $\theta \in {\Bbb R}$,
\[
R(\theta) = \left(
\begin{array}{cc}
\cos \theta & \sin \theta \\
- \sin \theta & \cos \theta
\end{array}
\right)
\]
and
\begin{equation}
\tilde{\omega} = exp ( i\, l\, \theta\, \tau' - i\, m\,
(\frac{\vec{v}^2}{2} \tau' + \vec{v} \cdot
\overrightarrow{R(\theta) u'}) - \frac{i k}{2} (\vec{v} \times
\overrightarrow{R(\theta) v'})) \label{7}
\end{equation}
Theorem II applies here as well. Therefore, it seems that only the case
$m=0$, $l \neq 0$, $k \neq 0$ has to be treated in the way indicated in
Ref.\cite{2}.

\section{The origin of cocycles}
We would like to understand the origin of nontrivial cocycles on Galilei
group $G$. It is the more interesting that the relativistic counterpart of $G$
-- the Poincare group $P$ -- does not admit nontrivial cocycles. On the
other hand, $G$ can be obtained from $P$ by contraction procedure. It is
therefore desirable to offer some interpretation for emergence of such
cocycles in nonrelativistic limit. The following general picture can be
given \cite{4}. Let $\omega(g,g')$ be any cocycle on $P$; write
\begin{equation}
\omega(g,g')=\exp i \xi(g,g') \label{8}
\end{equation}
Now, $\omega(g,g')$ is necessarily trivial, i.e. there exists a function
$\zeta$ on $P$ such that
\begin{equation}
\xi(g,g') = (\delta \zeta)(g,g') \equiv \zeta(g g') - \zeta(g) -
\zeta(g'). \label{9}
\end{equation}
The exponent $\xi(g,g')$ gives rise to a nontrivial cocycle in the
nonrelativistic limit $c\rightarrow \infty$ provided it survives the
contraction while $\zeta(g)$ does not (typically, it diverges as
$c\rightarrow \infty$). To make this pictures more concrete let us
describe in some detail the contraction procedure. First, we write the
element of Poincare group in matrix form
\begin{equation}
\{\Lambda, a\} \rightarrow \left[
\begin{array}{c|c}
\displaystyle \Lambda^{\mu}_{\nu} & \displaystyle a^{\mu}\\[1ex]
\hline
\displaystyle 0 & \displaystyle 1
\end{array}
\right] =
\left[
\begin{array}{c|c}
\displaystyle \delta^{\mu}_{\nu} & \displaystyle a^{\mu}\\[1ex]
\hline
\displaystyle 0 & \displaystyle 1
\end{array}
\right]
\left[
\begin{array}{c|c}
\displaystyle \Lambda^{\mu}_{\nu} & \displaystyle 0\\[1ex]
\hline
\displaystyle 0 & \displaystyle 1
\end{array}
\right]. \label{10}
\end{equation}
The Lorentz matrix $[\Lambda^{\mu}_{\nu}]$ is further decomposed into pure
boosts and rotations
\begin{equation}
\Lambda={\cal L}(\vec{v}) \cdot {\cal R} \label{11}
\end{equation}
where
\begin{equation}
{\cal L}(\vec{v}) =
\left[
\begin{array}{c|c}
\displaystyle \gamma & \displaystyle \frac{\gamma v_{k}}{c}\\[1ex]
\hline
\displaystyle \frac{\gamma v_{i}}{c} & \displaystyle \delta_{ik} + (\gamma
- 1) \frac{v_{i} v_{k}}{\vec{v}^2}
\end{array}
\right], \gamma \equiv \left(1-\frac{\vec{v}^2}{c^2}\right)^{-\frac{1}{2}}
\label{12a}
\end{equation}
\begin{equation}
{\cal R}=
\left[
\begin{array}{c|c}
\displaystyle 1 & \displaystyle 0\\[1ex]
\hline
\displaystyle 0 & \displaystyle R
\end{array}
\right], R R^{T}=R^{T} R=I. \label{12b}
\end{equation}
So, finally
\begin{equation}
\left[
\begin{array}{c|c}
\displaystyle \Lambda & \displaystyle a\\[1ex]
\hline
\displaystyle 0 & \displaystyle 1
\end{array}
\right]=
\left[
\begin{array}{c|c}
\displaystyle I & \displaystyle a\\[1ex]
\hline
\displaystyle 0 & \displaystyle 1
\end{array}
\right]
\left[
\begin{array}{c|c}
\displaystyle {\cal L}(\vec{v}) & \displaystyle 0\\[1ex]
\hline
\displaystyle 0 & \displaystyle 1
\end{array}
\right]
\left[
\begin{array}{c|c}
\displaystyle {\cal R} & \displaystyle 0\\[1ex]
\hline
\displaystyle 0 & \displaystyle 1
\end{array}
\right] \label{13}
\end{equation}
The contraction limit is now performed by multiplying eq.(\ref{13}) by $X$
from the right and $X^{-1}$ from the left, where
\begin{equation}
X=\left[
\begin{array}{c|c}
\displaystyle c & \displaystyle 0\\[1ex]
\hline
\displaystyle 0 & \displaystyle I
\end{array}
\right], \label{14}
\end{equation}
taking the limit $c \rightarrow \infty$ and identifying: $a^{0} \rightarrow
c \tau$, $\vec{a} \rightarrow \vec{u}$.

Now, one can easily explain the emergence of standard cocycle related to the
mass of particle. Take
\[\zeta(\{\Lambda,a\}) = c a^{0}\]
in eq.(\ref{9}). Due to the identification $a^{0} = c \tau$, $\zeta$
diverges as $c^2$ in the contraction limit. However,
\[
\zeta(\{\Lambda,a\}, \{\Lambda',a'\}) = c ( \Lambda^{0}_{\mu} {a'}^{\mu} +
a^{0} ) - c a^{0} - c {a'}^{0} =
c ( \Lambda^{0}_{\mu} - \delta^{0}_{\mu} ) {a'}^{\mu} =
\]
\begin{equation}
= c^2(\gamma-1)\gamma' + \gamma v_{i} R_{ik} u'_{k}
\stackrel{\scriptstyle c\rightarrow\infty}{\displaystyle \longrightarrow}
\frac{\vec{v}^2}{2} \cdot \tau' + \vec{v} \cdot
\overrightarrow{R u'}. \label{16}
\end{equation}
This explanation works both for three and four dimensions.\\
To account for the second cocycle (related to the parameter $k$) let us note
that, in the case of threedimensional space--time the rotation matrix
appearing in the decomposition (\ref{11})--(\ref{12b}) of the Lorentz
matrix is an element of $SO(2)$ and is therefore characterized by one angle
$\theta$. We put
\begin{equation}
\zeta(\{\Lambda,a\}) = c^2 \theta(\Lambda). \label{17}
\end{equation}
Actually, $\theta$ is multivalued no $P$ (while singlevalued on $\tilde{P}$)
but this plays no role in what follows. Now, from eq.(\ref{9}) we get
\begin{equation}
\xi(\Lambda,\Lambda') = c^2 ((\theta(\Lambda \cdot \Lambda') -
\theta(\Lambda) - \theta(\Lambda')). \label{18}
\end{equation}
But
\begin{equation}
\theta(\Lambda \cdot \Lambda') = \theta(\Lambda) + \theta(\Lambda') +
\delta\theta(\Lambda,\Lambda')  \label{19}
\end{equation}
where $\delta\theta(\Lambda,\Lambda')$ is $0(1/c^2)$ and is related to the
so called Thomas precession \cite{5}; its existence reflects the propety that
the composition of pure boosts is no longer a pure boost.\\

It follows from eqs.(\ref{17})--(\ref{19}) that $\xi$ survives the
$c\rightarrow\infty$ limit while $\zeta$ does not. To calculate
$\delta\theta$ we write
\[
\Lambda \cdot \Lambda' = ({\cal L}(\vec{v}){\cal R})({\cal L}(\vec{v'})
{\cal R}') =
{\cal L}(\vec{v})({\cal R}{\cal L}(\vec{v'}){\cal R}^{-1})({\cal R}{\cal
R'}) =
\]
\begin{equation}
= ({\cal L}(\vec{v}){\cal L}(\overrightarrow{R v'}))({\cal R}{\cal R}').
\label{20}
\end{equation}
The standard calculations (using eqs.(\ref{11}), (\ref{12a}), (\ref{12b}))
give
\begin{equation}
{\cal L}(\vec{v}){\cal L}(\overrightarrow{R v'}) = {\cal L}
(\overrightarrow{v''}){\cal R}(\delta\theta) \label{21}
\end{equation}
where the value of $\vec{v''}$ is there irrelevant while, in the limit
$c\rightarrow\infty$
\begin{equation}
\delta\theta = \frac{\vec{v} \times
\overrightarrow{R v'}}{2 c^2}. \label{22}
\end{equation}
By comparying eqs.(\ref{18}), (\ref{19}) and (\ref{22}) we get
\[
\xi = \frac{\vec{v} \times \overrightarrow{R v'}}{2}
\]
which gives the cocycle found previously.

\section{Bibliographical remarks}
We would like to conclude with the following bibligraphical remarks. The
central extensions of Lie algebra of threedimensional Galilei group were
found many years ago by Levy--Leblond \cite{6}. The corresponding cocycles
on Galilei group have been constructed by Grigore \cite{7}. In the same
paper Grigore has found the unitary projective representations of
$2+1$--dimensional Galilei group using the Mackey theory and exploiting the
trick (used again in Ref.\cite{2}) consisting in extending of Galilei group
with the help of two (three in the case of universal covering) central
charges. Grigore gave also a detailed discussion of projective
representations of $2+1$--dimensional Poincare group \cite{8}.


\newpage
\begin{thebibliography}{99}
\bibitem{1}Bose, S. K.: Comm. Math. Phys. \underline{169}, 385 (1995).
\bibitem{2}Bose, S. K.: Journ. Math. Phys. \underline{36}, 875 (1995).
\bibitem{3}Brihaye, Y., Gonera, C., Giller, S., Kosi\'{n}ski, P.:``Galilean
invariance in $2+1$ dimensions``, {\L}\'{o}d\'{z} University preprint (1994).
\bibitem{4}Saletan., E., J.: Journ. Math. Phys. \underline{2}, 1 (1961)\\
Aldaya, V., de Azcarraga, J., A.: Int. Journ. of Theor. Phys. \underline
{24}, 141 (1985).
\bibitem{5}M{\o}ller, C.: The Theory of relativity.
Oxford: Clarendon Press 1972.
\bibitem{6}L\'evy--Leblond., L--M.: Galilei group and Galilean
Invariance.\\
In: Group Theory and Its Applications. E. Loebl (ed.).\\ Academic Press
1971.
\bibitem{7}Grigore, D., R.: ``The projective unitary irreducible
representations of the Galilei group in $1+2$ dimensions``, preprint
IFA--FT--391--1993.
\bibitem{8}Grigore, D., R.: Journ. Math. Phys. \underline{34}, 4172 (1993).
\end{thebibliography}
\end{document}